\documentclass[10pt]{amsart}

\usepackage{amsmath, amssymb, amsthm, latexsym}

\theoremstyle{plain}
\newtheorem{theorem}{Theorem}[section]

\newtheorem{lemma}[theorem]{Lemma}
\newtheorem{corollary}[theorem]{Corollary}
\newtheorem{proposition}[theorem]{Proposition}
\newtheorem{definition}[theorem]{Definition}

\theoremstyle{remark}
\newtheorem{remark}[theorem]{Remark}

\input{xy}
\xyoption{all}

\newcommand{\cA}{\mathcal A}

\newcommand{\Tor}{{\rm Tor\,}}

\newcommand{\Pj}{{\rm Proj}}

\renewcommand{\dim}{{\rm dim}}
\renewcommand{\ker}{{\rm ker \,}}
\newcommand{\im}{{\rm im}}

\newcommand{\bm}{{\rm Bimod}}
\newcommand{\md}{{\rm Mod}}
\newcommand{\alg}{{\rm alg}}

\newcommand{\coker}{{\rm coker}\,}

\def\Cz{\mathbb{C}}
\def\Rz{\mathbb{R}}

\def\Nz{\mathbb{N}}
\def\Proj{{\rm Proj}}

\newcommand{\deq}{\stackrel{\rm def}{=}}     % equal by definition
\newcommand{\hot}{\mathbin{\overline{\otimes}}}
\newcommand{\hota}{\mathbin{\hat{\otimes}}}

\title{$L^2$-Invariants and rank metric}
\author{Andreas Thom}
\address{Andreas Thom, Mathematisches Institut der Universit\"at G\"ottingen,
Bunsenstr. 3-5, D-37073 G\"ottingen, Germany}
\email{thom@uni-math.gwdg.de}
\urladdr{http://www.uni-math.gwdg.de/thom}

\begin{document}

\begin{abstract}
We introduce a notion of rank completion for bi-modules over a finite tracial von Neumann algebra. We show that the functor of rank completion is exact and
that the category of complete modules is abelian with enough projective objects. This leads to interesting computations in the $L^2$-homology for
tracial algebras. As an application, we also give a new proof of a Theorem of Gaboriau on invariance of $L^2$-Betti numbers under orbit
equivalence.
\end{abstract}

\maketitle

\section{Preliminaries}

\subsection{Introduction}

The aim of this article is to unify approaches to several results in the theory of $L^2$-invariants of groups, see \cite{L, Gab}, 
and tracial algebras, see \cite{CS}. The new approach allows us to sharpen several results that were obtained in \cite{l2coh}. We also give a new
proof of D. Gaboriau's Theorem on invariance of $L^2$-Betti numbers under orbit equivalence.
In order to do so, we introduce the concept of rank metric and rank completion of bi-modules over a von finite tracial von Neumann algebra.

All von Neumann algebras in this article
have a separable pre-dual. Recall, a von Neumann algebra is called finite and tracial, if it comes with a fixed positive, faithful and normal trace. Every finite (i.e. Dedekind finite) von Neumann algebra admits such a trace, but we assume that a choice of a trace is fixed.

The rank is a natural measure of the size of the support of an element in a bi-module over a finite tracial von Neumann-algebra. 
The induced metric endows each bi-module with a topology, such that all bi-module maps are contractions. The main utility of completion with respect to the
rank metric is revealed by the observation that the functor of rank completion is exact and that the category of complete modules is abelian with
enough projective objects.

Employing the process of rank completion, we aim to proof two main results.
First of all, we will show that certain $L^2$-Betti number invariants of von Neumann algebras coincide with those for arbitrary weakly
dense sub-$C^*$-algebras. The particular case of the first $L^2$-Betti number was treated in \cite{l2coh}. The general result required a more
conceptual approach and is carried out in this article. The importance of this result was pointed out to the author by D. Shlyakhtenko. Indeed, according to A. Connes and D.
Shlyakhtenko, $K$-theoretic methods might be used to relate the $L^2$-Euler characteristic of a group $C^*$-algebra to the ordinary 
Euler characteristic of the group. This could finally lead to a computation of the $L^2$-Betti numbers for certain von Neumann algebras and would
resolve some longstanding conjectures, as for example the non-isomorphism conjecture for free group factors, see \cite{Voi}. 
However, a concrete implementation of this idea is not in reach and a lot preliminary work has still to be carried out.

Secondly, inspired by ideas of R. Sauer from \cite{sauer1}, we will give a new and self-contained 
proof of invariance of $L^2$-Betti numbers of groups under orbit equivalence. The idea here is very simple. We show that $L^2$-Betti number 
invariants cannot see the difference between an $L^{\infty}(X)$-algebra and its rank completion. Then, we observe the following: If free measure preserving actions 
of $\Gamma_1$ and $\Gamma_2$ on a probability space $X$ induce the same equivalence relation, then $L^{\infty}(X) \rtimes_{\alg} \Gamma_1$ and 
$L^{\infty}(X) \rtimes_{\alg} \Gamma_2$ have isomorphic rank completions as bi-modules (with respect to the diagonal left action) 
over $L^{\infty}(X)$. It remains to carry out several routine calculations in homological algebra.

\subsection{Dimension theory}

In his pioneering work, W. L\"uck was able to describe a lot of the analytic properties of the category of Hilbert-modules over a finite tracial von Neumann algebras $(M,\tau)$ in purely algebraic terms. This allowed to employ the maschinery of homological algebra in the study $L^2$-invariants and lead to substantial results and
a conceptional understanding from an algebraic point of view. One important ingredient in his work is
a dimension function which is defined for all $M$-modules, see \cite{L}. Due to several ring-theoretic properties of $M$, the natural dimension function
for projective modules has an extension to all modules and shares several convenient properies. In particular, it was shown in \cite{L}, that the sub-category
of zero-dimensional modules is a Serre sub-category, i.e. is closed under extensions. This implies that there is a $5$-Lemma for dimension isomorphisms. The following
lemma is immediate from this. (See \cite{W} for the necessary definitions.)

\begin{lemma} \label{deriv}
Let $(M,\tau)$ be a finite tracial von Neumann algebra.
Let $\cA$ be an abelian category with enough projective objects and let $F,G\colon \cA \to \md^M$ be right exact functors into the category $\md^M$ of $M$-modules. If there exists a natural transformation $h\colon F \to G$ which consists of dimension isomorphisms, then the induced natural transformations \[h_i\colon L_i(F) \to L_i(G)\] of left-derived functors consist of dimension isomorphisms too.
\end{lemma}

$L^2$-Betti numbers for certain group-actions on spaces 
were introduced by M. Atiyah in \cite{atiyah78}. The domain of definition was extended by J. Cheeger and M. Gromov in \cite{CG}. For references and
most of the main results, see \cite{L}. An important result of L\"uck was the following equality, which we take as a basis for our computations is Section \ref{gaboriau}: \begin{equation}\label{lueck} \beta^{(2)}_k(\Gamma) = \dim_{L\Gamma} \Tor_k^{\Cz\Gamma}(L\Gamma, \Cz).\end{equation}

The following observation concerning a characterization of zero-dimensional modules is due to R. Sauer, see \cite{sauer1}, and will be of major importance in the sequel.
\begin{theorem}[Sauer] \label{local}
Let $(M,\tau)$ be a finite tracial von Neumann algebra and let $L$ be a $M$-module. The following conditions are equivalent:
\begin{enumerate}
\item $L$ is zero dimensional.
\item $\forall \xi \in L,\, \forall \varepsilon >0, \, \exists p \in \Proj(M)\colon \quad \xi p = \xi \quad \mbox{and}\quad\tau(p) \leq \varepsilon.$
\end{enumerate}
\end{theorem}
The second condition is usually referred to as a local criterion of zero dimensionality. In the next section, we want to exploit this observation further and
study completions of bi-modules with respect to a certain metric that measures the size of the support.

\vspace{.2cm}

Let $a \in M$ be an arbitrary element in a finite tracial von Neumann algebra $(M,\tau)$. We denote by $s(a)$  is support projection and by $r(a)$ its range projection. Note that the equality $\tau(s(a)) = \tau(r(a))$ always holds. We denote by $\Pj(M)$ the set of projections of $M$. Note that $\Pj(M)$ is a complete, complemented modular
lattice. We denote the operations of meet, join and complement by $\wedge, \vee$ and $\perp$.

\section{Completion of bi-modules}
\subsection{Definition}
Let us denote the category of $M$-bi-modules by $\bm^M$. We will loosely identify $M$-bi-modules with $M \otimes M^o$-modules. 
In the sequel we regard $M \otimes M^o$ as a bi-module over $M$, acting by multiplication with $M \otimes M^o$ on the left.

\vspace{.2cm}
  
Let $(M,\tau)$ be a finite tracial von Neumann algebra and let $L$ be a $M$-bi-module. Associated to an element $\xi \in L$, there is a real-valued quantity that measures the size of the support.
Let us set
\[ [\xi] = \inf \left\{ \tau(p) + \tau(q)\colon p,q \in \Pj(M), p^\perp\xi q^\perp=0 \right\} \in [0,1]. \]
Obviously, for the bi-module $M$ and $x \in M$, we get that $[x]$ equals the trace of the support projection $s(x)$ of $x$. 
Indeed, if $p^{\perp}xq^{\perp} =0$, then $\tau(p^{\perp}) + \tau(s(x)) + \tau(q^{\perp}) \leq 2$ and thus
$$\tau(p) + \tau(q) \geq \tau(s(x)).$$ We conclude that $[x]\geq \tau(s(x))$. The reverse inequality is obvious.
\begin{lemma} \label{subadd}
Let $L$ be a $M$-bi-module and let $\xi_1, \xi_2 \in L$.
The in equality \[   [\xi_1 + \xi_2] \leq [\xi_1] + [\xi_2]\]
holds.
\end{lemma}
\begin{proof} Let $\varepsilon>0$ be arbitrary. We find projections $p_1,q_1,p_2,q_2$, such that
$\tau(p_i) + \tau(q_i) \leq [\xi_i] +\varepsilon$ and $p_i^{\perp}\xi_i q_i^{\perp}=0$, for $i=0,1$. Since
$(p_1^{\perp} \wedge p_2^{\perp})(\xi_1 + \xi_2)(q_1^{\perp} \wedge q_2^{\perp}) =0$ and, we get that
$$[\xi_1 + \xi_2] \leq \tau(p_1 \vee p_2) + \tau(q_1 \vee q_2) \leq \tau(p_1) + \tau(p_2) + \tau(q_1) + \tau(q_2) \leq
[\xi_1] + [\xi_2] + 2 \varepsilon.$$ Since $\varepsilon$ was arbitrary, the claim follows.
\end{proof}

Note, there is no reason to assume that $[\xi]=0 \Longrightarrow \xi=0$. Indeed, one can easily construct examples where this fails.

\begin{definition} Let $L$ be a $M$-bi-module. The quantity $d(\xi,\zeta) \deq [\xi -\zeta] \in \Rz$ defines a quasi-metric on $L$, which we
call \emph{rank metric}.
\end{definition}

\begin{lemma} \label{contraction}
Let $\phi\colon L \to L'$ be a homomorphism of $M$-bi-modules. 
\begin{enumerate}
\item The map $\phi$ is a contraction in the rank-metric.
\item Let $\varepsilon >0$ be arbitrary. If $\phi$ is surjective and $\xi' \in L'$, then there exists
$\xi \in L$, such that $\phi(\xi) = \xi'$ and $[\xi] \leq [\xi'] + \varepsilon$.
\end{enumerate}
\end{lemma}
\begin{proof}
(1) Let $\xi \in L$, $\varepsilon >0$ and $p,q \in \Pj(M)$, such that $\tau(p) + \tau(q) \leq [\xi] + \varepsilon$ and $p^{\perp}\xi q^{\perp}=0$.
Clearly, $p^{\perp} \phi(x) q^{\perp} =0$, and hence $[\phi(x)] \leq \tau(p) + \tau(q) \leq [\xi] + \varepsilon$. Since $\varepsilon$ was arbitrary, the assertion follows.

(2) Let $\xi' \in L'$. There exists $p,q \in \Pj(M)$, such that $\tau(p) + \tau(q) \leq [\xi'] + \varepsilon$ and $p^{\perp} \xi' q^{\perp}=0$.
Let $\xi'' \in L$ be any lift of $\xi'$ and set $\xi= \xi'' - p^{\perp} \xi'' q^{\perp}$. We easily see that $\phi(\xi) = \xi'$ and that
$p^{\perp}\xi q^{\perp} =0$. Hence, $\xi$ is a lift and $[\xi] \leq [\xi'] + \varepsilon$ as required.
\end{proof}

\begin{definition}
Let $L$ be a $M$ bi-module. The rank metric endows $L$ with a uniform structure.
\begin{enumerate}
\item We denote by $CS(L)$ the linear space of Cauchy sequences in $L$, by $ZS(L) \subset CS(L)$ the sub-space of sequences that converge to $0 \in L$.
Finally, we set $c(L) = CS(L)/ZS(L)$ and call it the \emph{completion} of $L$.
\item There is a natural map $L \to c(L)$ which sends an element to the constant sequence. The bi-module $L$ is called \emph{complete}, if it is an isomorphism. 
We denote by $\bm^M_c$ the full sub-category of complete $M \otimes M^o$-modules.
\end{enumerate}\end{definition}

\begin{lemma} Let $M$ be a finite tracial von Neumann algebra. 
\begin{enumerate}
\item The completion $M \hota M^o$ of $M \otimes M^o$ as a $M$-bimodule is a unital ring containing $M \otimes M^o$.
\item Let $L$ be a  $M$-bi-module. The completion is naturally a $M \hota M^o$-module and in particular a $M$-bi-module. 
\item The assignment $L \mapsto c(L)$ extends to a functor 
from the category of $M \otimes M^o$-modules to the category of $\bm^M_c$ of complete $M \otimes M^o$-modules.
\end{enumerate}
\end{lemma}
\begin{proof} Let $L$ be a $M$-bi-module. Let us first show that the $M \otimes M^o$-module structure extends to $c(L)$. 
Let $\xi \in M \otimes M^o$. We consider the map $\lambda_\xi: L \to L$ which is defined to be left-multiplication
by $\xi$. $\lambda_{\xi}$ is not a module-homomorphism but still to some extend compatible with the rank metric.
Let $\eta \in L, \varepsilon >0$ and $p,q \in \Pj(M)$ with $(p^{\perp} \otimes {q^{\perp o}})\eta=0$ and $\tau(p) + \tau(q) \leq [\eta] + \varepsilon$.

Specifying to $\xi = a \otimes b^o$, we get:
\[\lambda_{\xi}(\eta) = (a \otimes b^o)\eta = (a\otimes b^o)(1 \otimes 1 - p^{\perp} \otimes {q^{\perp o}}) \eta = (ap \otimes b^o) \eta + (ap^{\perp} \otimes (qb)^o) \eta . \]
We compute:
$(r(ap)^{\perp} \otimes 1)(ap \otimes b^o)\eta =0$ and hence $[(ap \otimes b^o)\eta] \leq \tau(r(ap)) = \tau(s(ap)) \leq \tau(p)$. Similarily,
$[(ap^{\perp} \otimes (qb)^o)\eta] \leq \tau(q)$ and hence $[(a \otimes b^o)\eta] \leq \tau(p) + \tau(q) \leq [\eta] + \varepsilon$. Again,
since $\varepsilon>0$ was arbitrary, we conclude $[\lambda_{a \otimes b^o}(\eta)] \leq [\eta]$.

If $\xi = a_1 \otimes b_1 + \dots + a_n \otimes b_n$, we get from Lemma \ref{subadd}, that
\[[L_{\xi}(\eta)] \leq n \cdot [\eta], \quad \forall \eta \in L.\]

We conclude that $\lambda_\xi$ is Lipschitz for all $\xi \in M \otimes M^o$. Hence, there is an extension $\lambda_{\xi}: CS(L)\to CS(L)$ which preserves
$ZS(L)$. Hence, there exists a bi-linear map $m': (M \otimes M^o) \times c(L) \to c(L)$ which defines a module structure that 
is compatible with the module structure on $L$.

It is clear that $[m(\xi,\eta)] \leq [\xi]$. Indeed, $(p^{\perp} \otimes q^{\perp o})\xi =0$ implies $(p^{\perp} \otimes q^{\perp o})m(\xi,\eta)=0$.
Hence $m'$ has a natural extension $$m: (M \hota M^o) \times c(L) \to c(L).$$
Obviously, if $L = M \otimes M^o$, then $m$ defines a multiplication that extends 
the multiplication on $M \otimes M^o$, i.e. the natural inclusion $M \otimes M^o \hookrightarrow M \hota M^o$ is 
a ring-homomorphism. This shows (1) and (2). Assertion (3) is obvious.
\end{proof}
\subsection{Completion is exact}

\begin{lemma} \label{exact}
The functor of completion is exact.
\end{lemma}
\begin{proof} 
Let $$  0 \to J \to L \stackrel{\pi}{\to} K \to 0$$ be an exact sequence of $M$-bi-modules. We have to show that
the induced sequence $$0 \to c(J) \to c(L) \to c(K) \to 0$$ is exact.

First, we consider the exactness at $c(K)$. Let $(\xi_n)_{n \in \Nz}$ be a Cauchy sequence in $K$. Without loss of generality, we can assume that
$[\xi_n - \xi_{n+1}] \leq 2^{-n}$. Lemma
\ref{contraction} implies that we can lift $(\xi_n)_{n \in \Nz}$ to a sequence $(\xi'_n)_{n \in \Nz} \subset L$ with
$[\xi'_n - \xi'_{n+1}] \leq 2^{1-n}$, hence a Cauchy sequence. This shows surjectivity of $c(\pi)$.

We consider now the exactness at $c(L)$. Obviously, $\im(c(J)) \subset \ker(c(\pi))$. 
Let $(\xi_n)_{n \in \Nz}$ be a Cauchy sequence in $L$ which maps to zero in $c(K)$. This says that
$(\pi(\xi_n))_{n \in \Nz}$ tends to zero. Again, by Lemma \ref{contraction}, we can lift $(\pi(\xi_n) )_{n \in \Nz}$ to a zero-sequence $(\xi'_n)_{n \in \Nz} \subset L$. Now, $(\xi_n - \xi'_n)_{n \in \Nz}$ defines a Cauchy sequence in $J$, that is equivalent to the sequence $(\xi_n)_{n \in \Nz}$ in the completion of $L$. Hence $\ker(c(\pi)) \subset \im(c(J))$ and the argument is finished.

The exactness at $c(J)$ is obvious since $J \subset L$ is a contraction in the rank metric by Lemma \ref{contraction}. This finishes the proof.
\end{proof}

\begin{theorem} \label{abel} Let $M$ be a finite tracial von Neumann algebra. Consider the category $\bm^M$ of $M$-bi-modules and 
the full sub-category $\bm^M_c$ of complete modules.
\begin{enumerate}
\item The completion functor $c\colon \bm^M \to \bm^M_c$ is left-adjoint to the forgetful functor from $\bm^M_c$ to $\bm^M$, i.e. whenever $K$ is complete:
$$\hom_{\bm^M}(c(L),K) = \hom_{\bm^M}(L,K).$$
\item The category $\bm^M_c$ is abelian and has enough projective objects.
\item The completion functor $c\colon \bm^M \to \bm^M_c$ preserves projective objects.
\item The kernel of the comparison map $L \to c(L)$ is $\{\xi \in L\colon [\xi]=0\}$.
\end{enumerate}
\end{theorem}
\begin{proof}
(1) If $K$ is complete, the natural map $K \to c(K)$ is an isomorphism, so that applying the functor $c$ defines a natural map
$$\hom_{\bm^M}(L,K) \to \hom_{\bm^M}(c(L),K).$$ A map in the inverse direction 
is provided by pre-composition with the map $L \to c(L)$. Assertion (1) follows easily by Lemma \ref{contraction}
since $\im(L)$ is dense in $c(L)$ and $\{\xi \in K, [\xi] =0\} = \{0\}.$

(2) If follows from Lemma \ref{exact}, that $\bm_c^M$ is abelian. Indeed, by exactness, kernels and co-kernels can be formed in $\bm^M$ and hence
all properties of those remain to be true in $\bm_c^M$. Let $L= \oplus_{\alpha} M \otimes M^o$ be a free $M \otimes M^o$-module. By (1), $c(L)$ is a projective object in $\bm_c^M$. If $K$ is complete and $\oplus_{\alpha} M \otimes M^o \stackrel{\pi}{\to} K$ is, using Lemma \ref{exact},  $c(\pi)$ is also a surjection onto $K$. I.e. there are
enough projective objects. (3) follows from (1). (4) is obvious.
\end{proof}
\subsection{Completion is dimension-preserving}
\begin{lemma} \label{map} Let $L$ be a $M$ bi-module. 
The natural map
\[{M \hot M^o} \otimes_{M \otimes M^o} L \to M \hot M^o\otimes_{M \otimes M^o} c(L)\]
is a dimension isomorphism.
\end{lemma}
\begin{proof} By Lemma \ref{exact}, it suffices to show that $$\ker(M \hot M^o \otimes_{M \otimes M^o} L \to M \hot M^o \otimes_{M \otimes M^o} c(L))$$ is zero-dimensional for all bi-modules $L$. 
Indeed $\coker(L \to c(L))$ has vanishing completion, and knowing the assertion for $\coker(L \to c(L))$ in place of $L$ implies that
$$\coker(M \hot M^o \otimes_{M \otimes M^o} L \to M \hot M^o \otimes_{M \otimes M^o} c(L))$$
has dimension zero.

\vspace{.2cm}

We want to apply the local critertion of Theorem \ref{local}. 
Let $\theta=\sum_{i=1}^n \eta_i \otimes \xi_i \in {M \hot M^o} \otimes_{M \otimes M^o} L$ be in the kernel.
This is to say that there exists some $l \in \Nz$ and zero-sequences $(\alpha_{i,k})_{k \in \Nz}$, for $1 \leq i \leq l$, such that
\[\sum_{i=1}^n \eta_i \otimes \xi_i = \sum_{i=1}^l \zeta_i \otimes \alpha_{i,k},\quad \forall k\in \Nz. \]
Indeed, the map factorizes through the split-injection
\[{M \hot M^o} \otimes_{M \otimes M^o} L \to M \hot M^o\otimes_{M \otimes M^o} CS(L)\]
and hence $\im({M \hot M^o} \otimes_{M \otimes M^o} ZS(L)) \subset M \hot M^o\otimes_{M \otimes M^o} CS(L)$ 
contains the image of
$$\ker(M \hot M^o \otimes_{M \otimes M^o} L \to M \hot M^o \otimes_{M \otimes M^o} c(L)).$$
Since $\alpha_{i,k} \to 0$, for all $1 \leq i \leq l$, for every $\varepsilon >0$
there exists $k$ big enough and projections $p_i,q_i \in \Pj(M)$, such that $(p_i^{\perp} \otimes q_i^{\perp o})\alpha_{i,k}  =0$ and 
$\tau(p_i) + \tau(q_i) \leq  \varepsilon/l$. 

Let $f_i$ be projections in $M \hot M^o$, such that $f_i^{\perp} \zeta_i = \zeta_i (p^{\perp}_i \otimes q^{\perp o})$. One can choose $f_i$ to satisfy
$\tau(f_i) \leq \tau(p) + \tau(q) \leq \varepsilon/l$, for all $1 \leq i \leq l$. 
We compute as follows:
\[ \bigwedge_{i=1}^l f_i^{\perp}\theta = \bigwedge_{i=1}^l f_i^{\perp} \left(\sum_{i=1}^l f_i^{\perp}\eta_i \otimes  \alpha_{i,k} \right)= 
\bigwedge_{i=1}^l f_i^{\perp} \left(\sum_{i=1}^l \eta_i \otimes  (p_i^{\perp} \otimes q_i^{\perp o}) \alpha_{i,k}\right) =0.\]
Thus $\tau(\vee_{i=1}^l f_i) \leq \varepsilon$ and and $(\vee_{i=1}^l f_i) \theta = \theta$.
Hence $$\ker\left({M \hot M^o} \otimes_{M \otimes M^o} L \to {M \hot M^o} \otimes_{M \otimes M^o} c(L)\right)$$ is zero-dimensional by Theorem \ref{local}.
\end{proof}

\begin{theorem} \label{ismo}
Let $\phi\colon L \to L'$ be a morphism of $M$-bi-modules. If $c(\phi)$ is an isomorphism, then 
$$\Tor_i^{M \otimes M^o}(M\hot M^o,L) \stackrel{\phi_*}{\to} \Tor_i^{M \otimes M^o}(M \hot M^o,L')$$ is a dimension isomorphism.
\end{theorem}
\begin{proof}
The exactness of $c\colon \bm^M \to \bm_c^M$ implies the following natural identification among left-derived functors:
$$L_i( M \hot M^o \otimes_{M \otimes M^o}?) \circ c = L_i(M \hot M^o \otimes_{M \otimes M^o}c(?)).$$
Indeed, this follows from the fact that $c\colon \bm^M \to \bm^M_c$ maps free modules in $\bm^M$ to projective objects 
in $\bm^M_c$. This implies the existence of a Grothendieck spectral sequence (see \cite[pp. 150]{W}) that yields the desired result.

Lemma \ref{map} together with Lemma \ref{deriv} implies the existence of a natural map
$$\Tor_i^{M \otimes M^o}(M\hot M^o,?) \to L_i( M \hot M^o \otimes_{M \otimes M^o} c(?))$$ which is a dimension isomorphism.
Combining the preceding two observations, we conclude that 
$$\dim_{M \hot M^o}\Tor_i^{M \otimes M^o}(M \hot M^o,L)=0, \quad \forall i \geq 0,$$ whenever $c(L) =0$. This implies the claim,
since $c(\ker(\phi))=0$ and $c(\coker(\phi))=0$ by Lemma \ref{exact}.
\end{proof}

\section{$L^2$-Betti numbers for tracial algebras}
\subsection{Preliminaries}
In \cite{CS}, A. Connes and D. Shlyakhtenko introduced a notion of $L^2$-homology and $L^2$-Betti numbers
for tracial algebras, compare also earlier work of W.L. Paschke in
\cite{Paschke}. The definition works well in a situation where the tracial algebra $(A,\tau)$ is contained in a finite von Neumann algebra $M$, to which the 
trace $\tau$ extends. More precisely, using the dimension function of W. L\"uck, see \cite{L}, they set:
\[\beta^{(2)}_k(A,\tau) = \dim_{M \hot M^o} \Tor_k^{A \otimes A^o}(M \hot M^o,A) \]

Here, $M \hot M^o$ denote the spatial tensor product of von Neumann algebras. We equipp $M \hot M^o$ with the trace $\tau \otimes \tau$.
Several results concerning these $L^2$-Betti numbers where obtained in \cite{CS} and \cite{l2sub,l2coh}. In particular, it was shown in \cite{CS} 
that $$\beta^{(2)}_k(\Cz\Gamma,\tau) = \beta^{(2)}_k(\Gamma),$$ where the right side denotes the $L^2$-Betti number of a group in 
the sense of Atiyah, see \cite{atiyah78} and Cheeger-Gromov, see \cite{CG}.

It is conjectured in \cite{CS} that $\beta^{(2)}_k(M,\tau)$ is an interesting invariant 
for the von Neumann algebra. Several related quantities where studied in
\cite{CS} as well. In particular,
\[\Delta_k(A,\tau) = \dim_{M \hot M^o} \Tor_k^{M \otimes M^o}(M \hot M^o,M\otimes_A M)\]
was studied for $k=1$.

\subsection{Pedersen's Theorem}
\begin{lemma} \label{dense}
Let $(M,\tau)$ be a finite tracial von Neumann algebra and let $A_1,A_2 \subset M$ be $*$-sub-algebras of $M$. If $A_1$ and $A_2$ have the same
closure with respect to the rank metric, then
\[\Delta_k(A_1,\tau) = \Delta_k(A_2,\tau),\quad \forall k\geq 0.\]
\end{lemma}
\begin{proof} Without loss of generality, $A_1 \subset A_2$. We show that $\pi \colon M\otimes_{A_1} M \to M \otimes_{A_2}M$ 
induces an isomorphism after completion with respect to the rank metric. The claim follows then from Theorem \ref{ismo}.

By Lemma \ref{exact}, it suffices to show that the kernel of $\pi$ has vanishing completion. An element $\xi$ in the kernel can be written as
\[\xi=\sum_{i=1}^n c_ia_i\otimes d_i - c_i \otimes a_id_i,\]
for some $c_i,d_i \in M$ and $a_i \in A_2$. Since $A_1$ is dense in $A_2$ there exists $a'_i \in A_1$ with $[a_i - a'_i] \leq \varepsilon/n$, for all
$1 \leq i \leq n$.

The following equality holds in $M \otimes_{A_1}M$:
\[\xi=\sum_{i=1}^n c_ia_i\otimes d_i - c_i \otimes a_id_i = \sum_{i=1}^n c_i(a_i - a'_i)\otimes d_i - c_i \otimes (a_i - a'_i)d_i,\]
Argueing as before, we see that there exists projections $p,q$ of trace less than $\varepsilon$, such that $p^\perp\xi q^\perp =0$. This implies that
$[\xi] \leq \varepsilon.$ Since $\varepsilon$ was arbitrary, we get that $[\xi] =0$ for all $\xi \in \ker(\pi \colon M
\otimes_{A_1} M \to M \otimes_{A_2} M)$ 
and thus $c(\ker(\pi))=0$. \end{proof}

The following result by G. Pedersen, see \cite[Thm. 2.7.3]{Ped}, is a deep result in the theory of operator algebras, which required a detailed analysis
of the precise position of a weakly dense $C^*$-algebra inside a von Neumann algebra. It is a generalization of a more classical theorem of 
Lusin in the commutative case.

\begin{theorem}[Pedersen]
Let $(M,\tau)$ be a finite tracial von Neumann algebra and let $A \subset M$ be a weakly dense sub-$C^*$-algebra. The algebra $A$ is dense in $M$
with respect to the rank metric.
\end{theorem}

\begin{corollary} \label{main} Let $(M,\tau)$ be a finite tracial von Neumann algebra and let $A \subset M$ be a weakly dense sub-$C^*$-algebra.
\[\Delta_k(A,\tau) = \Delta_k(M,\tau),\quad \forall k\geq 0.\]
\end{corollary}

\begin{remark}In \cite{l2coh}, it was shown that $\beta^{(2)}_1(A,\tau) = \beta^{(2)}_1(M,\tau)$, whenever $A$ is a weakly dense sub-$C^*$-algebra. In view of the factorization
\[\Tor_1^{A \otimes A^o}(M \hot M^o,A) \twoheadrightarrow \Tor_1^{M \otimes M^o}(M \hot M^o,M \otimes_A M) \to \Tor_1^{M \otimes M^o}(M \hot M^o,M ),\] the proof also shows that $\Delta_1(A,\tau) = \Delta_1(M,\tau)$ holds. 
Hence we can view Corollary \ref{main} as a generalization of this result from \cite{l2coh}.
\end{remark}

\section{Equivalence relations and Gaboriau's Theorem} \label{gaboriau}
\subsection{Equivalence relations and completion}

Most of the proofs in the section are parallel to proofs in Section $2$ and $3$ and hence we will give less detail and point to the relevant parts of Section $2$ and $3$.
Let $X$ be a standard Borel space and let $\mu$ be a probability measure on $X$. Given a discrete measurable equivalence relation (see \cite{FM1,FM2} for the 
necessary definitions)
\[R \subset X \times X,\]
we can form a \textit{relation ring} $R(X)$ as follows:
$$R(X) = \left\{\sum_{i=1}^n f_i \phi_i\colon f_i \in L^{\infty}(X), \phi_i \mbox{ local isomorphism from $R$} \right\} \subset L^{\infty}(R).$$
Here, $L^{\infty}(R)$ denote the generated von Neumann algebra, \see{FM1}. 
Note that $R(X)$ is a $L^{\infty}(X)$-bi-module with respect to the diagonal left action. (All $L^{\infty}(X)$-modules are bi-modules in this way.)
The following observation is the key to our results.

\begin{proposition} \label{key} Let $\Gamma$ be a discrete group and let $\rho\colon \Gamma \times X \to X$ be a measure preserving
free action of $G$ on $X$. We denote by $R_{\rho}$ the induced measurable equivalence relation on $X$.
The natural inclusion $\iota\colon L^{\infty}(X) \rtimes_{\alg} \Gamma \to R_{\rho}(X)$ induces an isomorphisms after completion.
\end{proposition}
\begin{proof}
According to the foundational work in \cite{FM1,FM2}, each local isomorphism $\phi$ which is implemented by the equivalence relation $R_{\rho}$ can be 
decomposed as a infinite sum of local isomorphism $\phi_i$, each of which is a cut-down of an isomorphism which is implemented by the action of a group element. 
The sizes of the supports of the cut-down local isomorphisms $\phi_i$ in this decomposition sum to the size of the 
support of $\phi$. It clearly implies, that $\phi$ can be approached by elements of 
$L^{\infty}(X) \rtimes \Gamma$ in rank metric. This finishes the proof.
\end{proof}

\begin{definition} Two group $\Gamma_1$ and $\Gamma_2$ are called \emph{orbit equivalent}, if there exists a probability space $X$ and free, measure preserving
actions of $\Gamma_1$ and $\Gamma_2$ on $X$ that induce the same equivalence relation.
\end{definition}

For an excellent survey on the properties of orbit equivalence and related notions, see \cite{Gabsurv}.

\begin{lemma} 
Let $R \subset X \times X$. We denote the completion of $R(X)$ by $\widehat{R}(X)$. $\widehat{R}(X)$ is a unital $L^{\infty}(X)$-algebra 
that contains $R(X)$ as a $L^{\infty}(X)$-sub-algebra.
\end{lemma}
\begin{proof} First of all, $\widehat{R}(X)$ is a $R(X)$-module. Indeed, left multiplication by $\sum_{i=1}^n f_i \phi_i$ is
easily seen to be Lipschitz with constant $n$. In particular, there exists a map $m'\colon R(X) \times \widehat{R}(X) \to \widehat{R}(X)$.

As before, we easily see that $m'$ has an extension to an associative and seperately continuous multiplication:
\[m\colon \widehat R(X) \times \widehat{R}(X) \to \widehat{R}(X). \] \end{proof}

\begin{lemma} \label{mod} 
Let $L$ be a $L^{\infty}(X) \rtimes_{\alg} \Gamma$-module. The completion of $L$, with 
respect to the diagonal $L^{\infty}(X)$-bi-module structure is naturally a $\widehat{R}(X)$-module and in particular a $R(X)$-module.
\end{lemma}
\begin{proof}
By Proposition \ref{key} $L^{\infty}(X) \rtimes_{\alg} \Gamma$ is dense in $R(X)$ and hence in $\widehat{R}(X)$.
Again, since $x = \sum_{i=1}^n f_i \phi_i$ acts with Lipschitz constant $n$ on $L$, the action extends to $c(L)$. Let $x_n$ be a Cauchy sequence in
$L^{\infty}(X) \rtimes\Gamma$ and $\xi \in c(L)$. The rank of $(x_n-x_m)\xi$ is less that the rank of $x_n-x_m$ and if $x_n \to 0$, then $x_n \xi \to 0$.
This finishes the proof.
\end{proof}

\subsection{Proof of Gaboriau's Theorem}
The proof of Gaboriau's Theorem which is presented in this section uses the technology of rank completion. It is very much inspired by
a proof of Gaboriau's Theorem given by R. Sauer in \cite{sauer1}. 
In his proof, the local criterion was a crucial ingredient to make the arguments work. We hope that the concept of rank completion will provide a good
understanding of why the Theorem is true.

\begin{lemma} \label{basics}
Let $\Gamma$ be a discrete group and let $\rho\colon \Gamma \times X \to X$ be a measurable and measure preserving action on a probability space $X$.
\begin{enumerate}
\item $\Cz\Gamma \subset L^{\infty}(X) \rtimes_{\alg} \Gamma$ is flat.
\item $L\Gamma \subset L^{\infty}(X) \rtimes \Gamma$ is flat and dimension preserving.
\end{enumerate} 
\end{lemma}
\begin{proof}
The first assertion is obvious. Indeed, $L^{\infty}(X) \rtimes_{\alg} \Gamma$ is a free $\Cz\Gamma$-module. The second assertion follows from the fact that
$L\Gamma$ is semi-hereditary, see \cite{sauer2}.
\end{proof}

\begin{theorem} \label{gab}
Let $\Gamma_1$ and $\Gamma_2$ are orbit equivalent groups, then
\[ \beta^{(2)}_k(\Gamma_1) = \beta^{(2)}_k(\Gamma_2),\quad \forall k\geq 0.\]
\end{theorem}
\begin{proof}
It suffices to write $\beta^{(2)}_k(\Gamma_1)$ entirely in terms of the equivalence relation it generates.
Using Lemma \ref{basics}, we rewrite:
\begin{eqnarray*}
(L^{\infty}(X) \rtimes \Gamma) \otimes_{L \Gamma}\Tor_k^{\Cz \Gamma}(L\Gamma,\Cz) 
&=&  \Tor_k^{\Cz \Gamma}(L^{\infty}(X) \rtimes \Gamma,\Cz)\\
&=& \Tor_k^{L^{\infty}(X) \rtimes_{\alg} \Gamma}(L^{\infty}(X) \rtimes \Gamma, L^{\infty}(X)).
\end{eqnarray*}
Here, the second equality follows since $$(L^{\infty}(X) \rtimes_{\alg} \Gamma) \otimes_{\Cz \Gamma} \Cz   = L^{\infty}(X)$$
as $L^{\infty}(X) \rtimes_{\alg} \Gamma$-module.
By Lemma \ref{basics} and Equation \ref{lueck}, we conclude that
\begin{equation}
\label{eq1}\beta_k^{(2)}(\Gamma) = \dim_{L^{\infty}(X) \rtimes \Gamma} \Tor_k^{L^{\infty}(X) \rtimes_{\alg} \Gamma}(L^{\infty}(X) \rtimes \Gamma,L^{\infty}(X)).
\end{equation}
There exists an exact functor which completes the category of $L^{\infty}(X) \rtimes_{\alg} \Gamma$-modules 
with respect to the diagonal left $L^{\infty}(X)$-bi-module structure. We have shown in Lemma \ref {mod} that the resulting full subcategory of those
$L^{\infty}(X) \rtimes_{\alg} \Gamma$-modules which are complete, is naturally a category of $R(X)$-modules. 
Let us denote the functor of completion by $$c\colon \md^{L^{\infty}(X) \rtimes_{\alg} \Gamma} \to \md_c^{R(X)}.$$

\begin{proposition} Let $L$ be a $L^{\infty}(X) \rtimes_{\alg} \Gamma$-module. The completion map induces an dimension isomorphism:
\[ (L^{\infty}(X) \rtimes \Gamma) \otimes_{L^{\infty}(X) \rtimes_{\alg} \Gamma} L \to (L^{\infty}(X) \rtimes \Gamma) \otimes_{R(X)} c(L).\]
\end{proposition}
\begin{proof}
The map can be factorized as 
\[ (L^{\infty}(X) \rtimes \Gamma) \otimes_{L^{\infty}(X) \rtimes_{\alg} \Gamma} L \to 
(L^{\infty}(X) \rtimes \Gamma) \otimes_{L^{\infty}(X) \rtimes_{\alg}\Gamma} c(L) \to
(L^{\infty}(X) \rtimes \Gamma) \otimes_{R(X)} c(L).\]
We show that each of the maps is a dimension isomorphism.
Let us start with the first one. Again, by exactness of $c$, it suffices to show that
$$\ker(L^{\infty}(X) \rtimes \Gamma) \otimes_{L^{\infty}(X) \rtimes_{\alg} \Gamma} L \to
L^{\infty}(X) \rtimes \Gamma) \otimes_{L^{\infty}(X) \rtimes_{\alg} \Gamma} c(L))$$ is zero-dimensional.
As in the proof of Lemma \ref{map}, an element $\theta$ in kernel is of the form:
\[\theta = \sum_{i=1}^l \zeta_i \otimes \alpha_{i,k},\quad \forall k \in \Nz,\]
for some zero sequences $(\alpha_{i,k})_{k \in \Nz} \subset L$. The proof proceeds as the proof of Lemma \ref{map}.

The second map can be seen to be a dimension isomorphism as follows. Clearly, the map is surjective and it remains to show
that the kernel is zero-dimensional. An element of the kernel is of the form:
\[\theta = \sum_{i=1}^n \xi_i \eta_i \otimes \zeta_i - \xi_i \otimes \eta_i \zeta_i, \]
for some $\xi_i \in L^{\infty}(X) \rtimes \Gamma, \eta_i \in R(X)$ and $\zeta_i \in L$. Approximating $\eta_i$ by elements in $L^{\infty}(X) \rtimes_{\alg} \Gamma$ we can assume (as in the proof of Lemma \ref{dense}) that $[\eta_i] \leq \varepsilon/(2n)$. The first summands are smaller than $\varepsilon/(2n)$, since support and range projection have the same trace. The second summand are also smaller than $\varepsilon/(2n)$ by the same argument 
and since projections in $L^{\infty}(X)$ can be moved through the tensor product. Hence $[\theta] \leq \varepsilon$. Since $\theta$ and $\varepsilon$ were arbitrary, we conclude by Theorem \ref{local} that $$\ker \left((L^{\infty}(X) \rtimes \Gamma) \otimes_{L^{\infty}(X) \rtimes_{\alg}\Gamma} c(L) \to
(L^{\infty}(X) \rtimes \Gamma) \otimes_{R(X)} c(L) \right)$$
is zero dimensional. This finishes the proof.
\end{proof}

To conclude the proof of Theorem \ref{gab}, we note that by Lemma \ref{deriv} we get an induced map
$$\Tor_i^{L^{\infty}(X) \rtimes_{\alg} \Gamma}(L^{\infty}(X) \rtimes \Gamma,?) \to
L_i((L^{\infty}(X) \rtimes \Gamma) \otimes_{R(X)} c(?) )$$ which is a dimension isomorphism. 
The right hand side applied to $L^{\infty}(X)$ depends only on the generated equivalence relation. Indeed, as in the proof of Theorem \ref{ismo}, a Grothendieck
spectral sequence shows $$L_i((L^{\infty}(X) \rtimes \Gamma) \otimes_{R(X)} c(?) ) = L_i((L^{\infty}(X) \rtimes \Gamma) \otimes_{R(X)} ? ) \circ c.$$
Here, we use implicitly that the category of complete $R(X)$-modules is abelian with enough projective objects. The proof of this fact can be taken verbatim
from the proof of Theorem \ref{abel} and the adjointness relations:
\[\hom_{L^{\infty}(X) \rtimes \Gamma}(L,K) = \hom_{L^{\infty}(X) \rtimes \Gamma}(c(L),K) = \hom_{R(X)}(c(L),K).\]
The projective objects are completions of free $R(X)$-modules. 
\end{proof}

\begin{remark} There is a second major result of Gaboriau's on proportionality of $L^2$-Betti numbers for weakly orbit equivalent groups, see \cite{Gabsurv}. 
Sauer has shown in \cite{sauer1}, how homological methods and properties of the dimension function allow to deduce this result. The same arguments
apply to our setting.
\end{remark}
\bibliographystyle{alpha}
\bibliography{l2.bib}

\end{document}